\documentclass[twoside]{amsart}

\usepackage{latexsym,amssymb}

\def\demo{{\it Proof. }}
\def\sqr#1#2{{\vcenter{\hrule height.#2pt
        \hbox{\vrule width.#2pt height#1pt \kern#1pt
                \vrule width.#2pt}
        \hrule height.#2pt}}}
\def\square{\mathchoice\sqr64\sqr64\sqr{4}3\sqr{3}3}
\def\QED{\hfill$\square$}

\newtheorem{Theorem}{\indent \sc Theorem}[section]
\newtheorem{Lemma}[Theorem]{\indent \sc Lemma}
\newtheorem{Remark}[Theorem]{\indent \sc Remark}
\newtheorem{Example}[Theorem]{\indent \sc Example}

\begin{document}

\title[Strongly Cohen--Macaulay ideals]
{\Large Strongly Cohen--Macaulay ideals \\
of small second analytic deviation}

\author[A. Corso and C. Polini]
{Alberto Corso \and Claudia Polini}

\thanks{
AMS 2000 {\it Mathematics Subject Classification}. Primary 13H10;
Secondary 13C40,  13D02, 13D25. \newline\indent
The NSF, under grant DMS-9970344, has partially supported the
research of the second author and has therefore her heartfelt thanks. }

\address{Department of Mathematics, Michigan State University, E.
Lansing, Michigan 48824}
\email{corso@math.msu.edu}

\address{Department of Mathematics, University of Oregon, Eugene,
Oregon 97403}
\email{polini@math.uoregon.edu}

\dedicatory{To Professor William Heinzer  on the occasion of his
sixtieth birthday}

\begin{abstract}
We characterize the strongly Cohen--Macaulay ideals of second
analytic deviation one in terms of depth properties of the powers of
the ideal in the `standard range.' This provides an explanation of the
behaviour of certain ideals that have appeared in the literature.
\end{abstract}

\maketitle

\section{Introduction}

Let $R$ be a local Gorenstein ring. Let $I$ be an $R$-ideal of grade
$g$ and minimally generated by $n$ elements.
Let $H_{*}(I)$ denote the homology modules of the ordinary Koszul complex
${\mathbb K}_{*}$ built on a minimal generating set of $I$.
The ideal $I$ is said to be {\it strongly Cohen--Macaulay} if all the
$($non-zero$)$ homology modules $H_i(I)$ are Cohen--Macaulay.
The notion of strong Cohen--Macaulayness appeared as a generalization and
correction by C. Huneke \cite{Hscm} of a result of M. Artin and
M. Nagata \cite{AN} on residual Cohen--Macaulayness.

An extreme case of this property is given by a complete intersection.
A Cohen--Macaulay almost complete intersection is also strongly
Cohen--Macaulay. In \cite{AH}, L. Avramov
and J. Herzog settled the next case: Any Cohen--Macaulay ideal that can
be generated by $g+2$ elements is strongly
Cohen--Macaulay.
Other distinguished classes of strongly Cohen--Macaulay ideals include
perfect ideals of codimension two and perfect Gorenstein ideals of
codimension three. More generally, ideals in the linkage class of a complete
intersection, a situation that was dubbed {\it licci ideals}, are strongly
Cohen--Macaulay \cite{Hlkhi}.

The property of being strongly Cohen--Macaulay plays a crucial role in the
study of {\it blowup algebras}.
We recall $I$ satisfies property $G_{s}$, where $s$ is an integer, if
$\mu(I_{\mathfrak p}) \leq {\rm dim}\, R_{\mathfrak p}$ for any prime ideal
${\mathfrak p}$ containing $I$ with ${\rm dim}\, R_{\mathfrak p} \leq s-1$.
We say  $I$ satisfies $G_{\infty}$ if $I$ satisfies $G_s$ for all $s$.
A classical result of \cite{HSV}, later refined in \cite{HVV}, states that:
If $I$ is a strongly Cohen--Macaulay ideal satisfying
$G_{\infty}$ then the symmetric algebra ${\rm Sym}(I)$ and the Rees algebra
$R[It]$ of $I$ are isomorphic. Similarly, ${\rm Sym}(I/I^2)$ is isomorphic
to the associated graded ring ${\rm gr}_I(R)$ of $I$ and, in both cases,
these algebras are all Cohen--Macaulay. In this situation, the ideal $I$ is
said to be of {\it linear type}.

In this paper we consider ideals $I$ with property $G_{\ell}$,
where $\ell$ is the analytic spread of $I$, minimally generated by
$n=\ell+1$ elements, with $\ell \geq g+1$ and satisfying ${\rm
depth}\, R/I^j \geq {\rm dim}\, R/I-j+1$ for $1 \leq j \leq
\ell-g$. Ideals for which $n=\ell+1$ are usually referred to as
ideals of {\it second analytic deviation} one \cite{UV}.
In presence of property $G_{\ell}$, it is also well known that the
above depth conditions are satisfied if the ideal is strongly
Cohen--Macaulay \cite{U}. Ideals of second analytic
deviation one have already been extensively studied in the
literature $($see for example \cite{AH2,AHH,J,PU2,UV}$)$ and to
give a flavor of their relevance we quote a result of \cite{PU2},
which states that: The only perfect Gorenstein $($not of linear type$)$
ideals of codimension three with $G_{\ell}$ and with Cohen--Macaulay
Rees algebra are exactly the ideals of second analytic deviation one.

The goal of this paper is to establish the strongly Cohen--Macaulay
property of the whole class of such ideals $($see
Theorem~\ref{scm}$)$, thereby providing a better understanding of recent
results from, for example, \cite{J} and \cite{PU2}. At the same time,
examples are given to illustrate that the theorem breaks down if
either the second analytic deviation requirement or the depth
conditions are relaxed.

\section{The main result}

In order to prove the main result we first need to study the degrees of
the defining equations of the Rees algebra of our ideals.
This will be done in Lemma~\ref{syzygetic}, after we have established some
definitions and facts.

There are two natural exact sequences that are useful to represent
externally the Rees algebra of an ideal $I$. The first one
\[
0 \rightarrow {\mathcal A} \longrightarrow {\rm Sym}(I) =
\bigoplus_{j=0}^{\infty} S_j(I) \longrightarrow R[It] = \bigoplus_{j=0}^{\infty}
I^jt^j \rightarrow 0
\]
relates the symmetric algebra ${\rm Sym}(I)$ and the Rees algebra $R[It]$ of
$I$. Since those two objects are graded, we can also write
${\mathcal A}$ in terms of its graded components, namely
${\mathcal A} = {\mathcal A}_{\,2} + {\mathcal A}_{\,3} + \cdots$.
Another way of representing the Rees algebra of $I$ is by means of the
following mapping
\[
0 \rightarrow {\mathcal Q} \longrightarrow R[T_1, \ldots, T_n]
\longrightarrow R[It] \rightarrow 0,
\]
where $I = (a_1, \ldots, a_n)$ and $T_j$ is mapped to $a_jt$.
${\mathcal Q}$ is the ideal generated by the equations of the Rees algebra.
Again write ${\mathcal Q}$ in terms of its graded components, namely
${\mathcal Q} = {\mathcal Q}_{\, 1} + {\mathcal Q}_{\, 2} + \cdots$.
The above short exact sequences are related in the
following manner
\[
{\rm Sym}(I) \simeq B/B{\mathcal Q}_{\,1} \qquad
{\mathcal A}_{\,j} = {\mathcal Q}_{\,j}/B_{j-1}{\mathcal Q}_{\,1} \quad
{\rm for \quad } j \geq 2,
\]
where $B=R[T_1, \ldots, T_n]$.
Therefore ${\mathcal A}_{\,2} = \ldots = {\mathcal A}_{\,r} = 0$ is equivalent
to saying that $S_j(I) \simeq I^j$ for $1 \leq j \leq r$ and also that
the equations of the Rees algebra up to degree $r$ are generated by the
linear ones.

On the other hand, it is also useful to study internally the Rees algebra
of $I$, that is by means of reductions of $I$.
A subideal $J$ of $I$ is said to be a reduction of $I$ if
the inclusion of Rees algebras $R[Jt] \hookrightarrow R[It]$ is module
finite \cite{NR}.
The philosophy underlying the use of reductions is that it is reasonable
to expect to recover some of the properties of $R[It]$ from those
of $R[Jt]$. It is even more reasonable to expect better results on
$R[It]$ when $J$ has an amenable structure. Minimal
reductions are reductions minimal with respect to containment.
If the residue field is infinite, their minimal number of generators does
not depend on the minimal reduction of the ideal. This number $\ell$
is called analytic spread of $I$. It equals the dimension of the fiber
cone of $I$ and it is always greater than or equal to the height of $I$.

Finally, a proper $R$-ideal $K$ is called
an $s$-residual intersection of $I$, where $s$ is an integer,
if there exists an $s$-generated ideal ${\mathfrak a} \subset I$ so
that $K= {\mathfrak a} \colon I$ and the height of $K$ is at least $s$.
If in addition the height of $I + K$ is at least $s+1$ we say that $K$
is a geometric $s$-residual intersection of $I$.
The next lemma summarizes all the properties of residual intersections
that will be used in the sequel of the paper. We remark that these results
are valid in a much broader situation, for which we refer to
\cite{U,JU,CP}. The way we stated them has been tailored to fit
our present needs.

\begin{Lemma}\label{properties}
Let $R$ be a local Gorenstein ring, let $K = {\mathfrak a} \colon I$ be
an $s$-residual intersection of $I$ and assume that $I$ satisfies $G_s$.
\begin{itemize}
\item[$(${\it a}$)$]
There exists a generating sequence $a_1, \ldots, a_s$ of ${\mathfrak a}$
such that $K_i = {\mathfrak a}_i \colon I$, where ${\mathfrak a}_i =
(a_1, \ldots, a_i)$, are geometric $i$-residual intersections of $I$
for $ i \leq s-1$.
\end{itemize}
Assume further that ${\rm depth}\, R/I^j \geq {\rm dim}\, R/I -j+1$ for
$1 \leq j \leq s-g$.
\begin{itemize}
\item[$(${\it b}$)$]
Any sequence $a_1, \ldots, a_s$ as in $(${\it a}$)$ forms an unconditioned
$d$-sequence.

\item[$(${\it c}$)$]
$K_i = {\mathfrak a}_i \colon (a_{i+1})$ for $i \leq s-1$.

\item[$(${\it d}$)$]
${\rm depth}\, R/{\mathfrak a}_i = d-i$ for $i \leq s$.

\item[$(${\it e}$)$]
$K_i \cap I^t = {\mathfrak a}_iI^{t-1}$ for $i \leq s-1$ and $1
\leq t \leq s-g$.

\item[$(${\it f}\/$)$]
${\mathfrak a}I^t \colon I^t = {\mathfrak a}$ for $ 1 \leq t \leq
s-g$.
\end{itemize}
\end{Lemma}
\demo See \cite[1.6, 1.7, 1.8 and 2.9]{U} for a proof of $(${\it
a}$)$-$(${\it d}$)$. See \cite[2.5]{JU} together with \cite[2.9]{U}
for a proof of $(${\it e}$)$.
See \cite[2.4]{CP} for a proof of $(${\it f}\/$)$. \QED

\medskip

The typical situation in which we apply Lemma~\ref{properties} is
described next $($see \cite[1.11 and 2.9]{U} or \cite[2.7]{JU}
together with \cite[2.9]{U}$)$.

\begin{Remark}\label{remark}
{\rm Let $R$ be a local Gorenstein ring with infinite residue field.
Let $I$ be an $R$-ideal with analytic spread $\ell$ satisfying $G_{\ell}$
and ${\rm depth}\, R/I^j \geq {\rm dim}\, R/I -j+1$ for
$1 \leq j \leq \ell-g$. Then for any minimal reduction $J$ of $I$
the ideal $K = J \colon I$ is an $\ell$-residual intersection of $I$. }
\end{Remark}

Lemma~\ref{syzygetic} below, which contains our main computation,
is a very interesting result in its own right as it shows that, in
our settings, the Rees algebra of $I$ has no equations in low
degrees. This result simplifies earlier work of Johnson, Polini
and Ulrich, who assumed the condition ``$S_j(I) \simeq I^j$ for
every $1 \leq j \leq \ell-g+1$'' in several results of \cite{J}
and \cite{PU2}. Our proof was in part suggested by the one of
\cite[2.8]{H}, which in turn uses {\it ad hoc} techniques for
second analytic deviation one ideals used in earlier works
$($see \cite{AH2}, for example$)$.

\begin{Lemma}\label{syzygetic}
Let $R$ be a local Gorenstein ring with infinite residue field.
Let $I$ be an $R$-ideal of grade $g$ and satisfying $G_{\ell}$, where $\ell$
is the analytic spread of $I$, and ${\rm depth}\, R/I^j \geq {\rm dim}\, R/I
-j+1$ with $1 \leq j \leq \ell-g$. If $I$ has second analytic deviation
one then $S_j(I) \simeq I^j$ for every $1 \leq j \leq \ell-g+1$.
\end{Lemma}
\demo
Write $I = (a_1, \ldots, a_{\ell+1})$ and assume that $J = (a_1, \ldots,
a_{\ell})$ is a minimal reduction of $I$. Also, assume that the $a_i$
have been chosen as in Lemma~\ref{properties}.
Finally, we may assume that $\ell \geq g+1$.

With the setup and notations established before, we first show by
descending induction on $q$, where $1 \leq q \leq
\ell$, that
\[
{\mathcal Q}_{\,j} \subset {\mathcal Q}_{\,j} \cap
(T_1, \ldots , T_{\ell+1})^{j-1}(T_1, \ldots, T_{q-1}) +
(T_q, \ldots, T_{\ell}){\mathcal Q}_{\,j-1} +
T_{\ell+1}^{j-1}{\mathcal Q}_{\,1},
\]
for all $2 \leq j \leq \ell-g+1$.

\noindent \underline{\it Base of induction$:$} we show the statement for $q =
\ell$. Let $H$ be an homogeneous polynomial of degree $j$, $H \in
R[T_1, \ldots, T_{\ell+1}]$, such that $H(a_1, \ldots, a_{\ell+1})
= 0$. For some $0 \leq e \leq j$ we can write
\begin{equation}\label{eq1}
H = r(T_1, \ldots, T_{\ell})T_{\ell+1}^e + p(T_1, \ldots, T_{\ell+1}),
\end{equation}
where $r(T_1, \ldots, T_{\ell})$ is a homogeneous polynomial of degree $j-e$
in $T_1, \ldots, T_{\ell}$ and $p(T_1, \ldots, T_{\ell+1})$ is a
homogeneous polynomial of degree $j$ in $T_1, \ldots, T_{\ell+1}$
having degree at most $e-1$ in $T_{\ell+1}$.

\noindent \underline{\it Case {\rm 1:}} $e = j$.
Evaluate equation $($\ref{eq1}$)$ by substituting the $T$'s with the
corresponding $a$'s. We get
\[
r a_{\ell+1}^j = p(a_1, \ldots, a_{\ell+1})
\in J I^{j-1}
\]
with $r \in R$. Thus $r I^j \subset JI^{j-1}$ with $1 \leq j-1
\leq \ell - g$. Hence $r I \subset JI^{j-1} \colon I^{j-1} = J$ by
Lemma~\ref{properties}$(${\it f}\/$)$. We can find a linear
polynomial of the form
\[
L = L(T_1, \ldots, T_{\ell+1}) = r T_{\ell+1} + L'(T_1, \ldots, T_{\ell}),
\]
with $L'(T_1, \ldots, T_{\ell})$ linear in $T_1, \ldots,
T_{\ell}$, such that $L(a_1, \ldots, a_{\ell+1}) = 0$.
Consider the polynomial $P = P(T_1, \ldots, T_{\ell+1}) =
H - T_{\ell+1}^{j-1}L$.
We have that $P(a_1, \ldots, a_{\ell+1}) = 0$ and moreover $P$ can be
written as
\begin{equation}\label{eq2}
P = T_1 P_1 + \ldots + T_{\ell}P_{\ell}
\end{equation}
with $P_k \in  (T_1, \ldots, T_{\ell+1})^{j-1}$. Set $p_k =
P_k(a_1, \ldots, a_{\ell+1})$. Evaluate equation $($\ref{eq2}$)$
by substituting the $T$'s with the corresponding $a$'s. We obtain,
using Lemma~\ref{properties}$(${\it c}$)$ and $(${\it e}$)$,
\[
p_{\ell} \in ((a_1, \ldots, a_{\ell-1}) \colon a_{\ell}) \cap I^{j-1} =
(a_1, \ldots, a_{\ell-1}) I^{j-2}.
\]
Hence we can find another polynomial $P_{\ell}^{*} =
P_{\ell}^{*}(T_1, \ldots, T_{\ell+1}) = T_1 H_1 +
\ldots + T_{\ell-1}H_{\ell-1}$, for some $H_k \in  (T_1, \ldots,
T_{\ell+1})^{j-2}$, such that $P_{\ell}^{*}(a_1, \ldots, a_{\ell+1}) =
p_{\ell}$. Thus $P_{\ell}-P_{\ell}^{*} = \widetilde{P}_{\ell}
\in {\mathcal Q}_{\,j-1}$. Hence $P$ can be rewritten as
\begin{eqnarray*}
P & = & T_1 P_1 + \ldots + T_{\ell-1}P_{\ell-1} +
T_{\ell}(P_{\ell}^{*} + \widetilde{P}_{\ell}) \\
& = & T_1 P_1 + \ldots + T_{\ell-1}P_{\ell-1} +
T_{\ell}(T_1 H_1 +
\ldots + T_{\ell-1}H_{\ell-1}+ \widetilde{P}_{\ell}) \\
& = & T_1 (P_1+T_{\ell}H_1) + \ldots + T_{\ell-1}(P_{\ell-1}
+T_{\ell}H_{\ell-1}) + T_{\ell}\widetilde{P}_{\ell} \\
& = & T_1 P_1' + \ldots + T_{\ell-1}P_{\ell-1}' + T_{\ell}\widetilde{P}_{\ell}
\end{eqnarray*}
for some $P_k' \in (T_1, \ldots, T_{\ell+1})^{j-1}$. In conclusion we have
that $H = T_1 P_1' + \ldots + T_{\ell-1}P_{\ell-1}' +
T_{\ell}\widetilde{P}_{\ell} + T_{\ell+1}^{j-1} L$, which shows that
\[
H \in {\mathcal Q}_{\,j} \cap
(T_1, \ldots , T_{\ell+1})^{j-1}(T_1, \ldots, T_{\ell-1}) +
T_{\ell} {\mathcal Q}_{\,j-1} +
T_{\ell+1}^{j-1}{\mathcal Q}_{\,1}.
\]

\noindent \underline{\it Case {\rm 2:}} $e \leq j-1$. In this case we have that
$H \in Q_{\,j} \cap (T_1, \ldots, T_{\ell+1})^{j-1}(T_1, \ldots,
T_{\ell})$. Hence we can write $H = T_1 H_1 + \ldots + T_{\ell}H_{\ell}$
with $H_k \in  (T_1, \ldots, T_{\ell+1})^{j-1}$. Use the same argument
we used for $P$ in the previous case and conclude that
\[
H = T_1 H_1' + \ldots + T_{\ell-1}H_{\ell-1}' + T_{\ell}\widetilde{H}_{\ell}
\]
for some $H_k' \in (T_1, \ldots, T_{\ell+1})^{j-1}$ and $\widetilde{H}_{\ell}
\in Q_{\,j-1}$.

\noindent \underline{\it Inductive step$:$} Suppose that we have
\[
{\mathcal Q}_{\,j} \subset {\mathcal Q}_{\,j} \cap
(T_1, \ldots , T_{\ell+1})^{j-1}(T_1, \ldots, T_{q-1}) +
(T_q, \ldots, T_{\ell}){\mathcal Q}_{\,j-1} +
T_{\ell+1}^{j-1}{\mathcal Q}_{\,1}.
\]
We want to show that
\[
{\mathcal Q}_{\,j} \subset {\mathcal Q}_{\,j} \cap
(T_1, \ldots , T_{\ell+1})^{j-1}(T_1, \ldots, T_{q-2}) +
(T_{q-1}, \ldots, T_{\ell}){\mathcal Q}_{\,j-1} +
T_{\ell+1}^{j-1}{\mathcal Q}_{\,1}.
\]
Clearly, we only need to focus on terms $H \in {\mathcal Q}_{\,j} \cap
(T_1, \ldots , T_{\ell+1})^{j-1}(T_1, \ldots, T_{q-1})$.
Write $H = T_1 H_1 + \ldots + T_{q-1}H_{q-1}$
with $H_k \in  (T_1, \ldots, T_{\ell+1})^{j-1}$. Use the same argument
we used for $P$ in the previous case and conclude that
\[
H = T_1 H_1' + \ldots + T_{q-2}H_{q-2}' + T_{q-1}\widetilde{H}_{q-1}
\]
for some $H_k' \in (T_1, \ldots, T_{\ell+1})^{j-1}$ and $\widetilde{H}_{q-1}
\in Q_{\,j-1}$.

\noindent \underline{\it Conclusion$:$} for any $2 \leq j \leq \ell-g+1$ we
have therefore that
\[
{\mathcal Q}_{\,j} \subset (T_1, \ldots, T_{\ell}) {\mathcal Q}_{\,j-1}+
T_{\ell+1}^{j-1}{\mathcal Q}_{\,1} \subset (T_1, \ldots, T_{\ell+1})
{\mathcal Q}_{\,j-1}.
\]
By iteration we obtain ${\mathcal Q}_{\,j} \subset (T_1, \ldots,
T_{\ell+1})^{j-1} {\mathcal Q}_{\,1}$, as desired. \QED

\bigskip

We are ready to prove our main result: the key point is to show
that the apparently weaker condition stated in $(${\it b}$)$ is indeed
equivalent to the one in $(${\it a}$)$.

\begin{Theorem}\label{scm}
Let $R$ be a local Gorenstein ring with infinite residue field.
Let $I$ be an $R$-ideal of grade $g$ and satisfying $G_{\ell}$, where $\ell$
is the analytic spread of $I$. Suppose further that $I$ has second
analytic deviation one and $\ell \geq g+1$.
Then the following conditions are equivalent
\begin{itemize}
\item[$(${\it a}$)$]
$I$ is a strongly Cohen--Macaulay ideal;

\item[$(${\it b}$)$]
${\rm depth}\, R/I^j \geq {\rm dim}\, R/I -j+1$ for $1 \leq j \leq \ell-g$.
\end{itemize}
\end{Theorem}
\demo That $(${\it a}$)$ implies $(${\it b}$)$ is a general result
proved in \cite[2.10]{U}. We are concerned with the proof of the
converse.

By adjoining a set of indeterminates to $R$ and $I$ and after localizing, we
may assume the grade $g$ of $I$ is at least $2$. This leaves all assumptions
preserved. The Koszul homology is unchanged: the net effect of adding
variables is that we have a new Koszul complex whose exact tail is longer.
Also, the lower bounds on the depths of the powers of the ideal are still
satisfied.

Write $I=(a_1, \ldots, a_{\ell}, a_{\ell+1})$ and assume that $J =
(a_1, \ldots, a_{\ell})$. Also, assume that the $a_i$ have been
chosen as in Lemma~\ref{properties}. By
Lemma~\ref{properties}$(${\it b}$)$ $J$ is generated by a
$d$-sequence, hence by a proper sequence $($see \cite[2.3.1$(${\it
c}$)$]{red}$)$. We now claim that the ideal $I$ is generated by a
proper sequence. For that, it is enough to show that $I$
annihilites the Koszul homology modules of $J$. Notice that
Lemma~\ref{properties}$(${\it d}$)$ and \cite[3.7]{HVV} imply that
the ideal $J$ satisfies sliding depth, that is ${\rm depth}\,
H_j(J) \geq d-\ell+j$ for $0 \leq j \leq \ell-g$ and $d={\rm
dim}\, R$. Now, since $I H_j(J)$ is a submodule of $H_j(J)$ and
${\rm depth}\, H_j(J) \geq d-\ell+j$ for $1 \leq j \leq \ell-g$,
it suffices to check $IH_j(J)=0$ at prime ideals ${\mathfrak p}$
of codimension at most $\ell-j \leq \ell-1$. But in this range
$I_{\mathfrak p} = J_{\mathfrak p}$.

To conclude that $H_j(I)$ is Cohen--Macaulay for $0 \leq j \leq
\ell-g+1 = n-g$, we need to estimate the depths of the modules of
cycles and boundaries $Z_j$ and $B_j$ of the associated Koszul
complex ${\mathbb K}_{*}$ of $I$. Since $I$ is generated by a
proper sequence, it follows from \cite[12.5]{HSV} that the
approximation complex ${\mathcal Z}_{*}$ on $a_1, \ldots,
a_{\ell+1}$ is acyclic with $H_0({\mathcal Z}_{*}) \simeq {\rm
Sym}(I)$. Let $F=R^n$ and let
\[
0 \rightarrow Z_j \otimes S_0(F)
\rightarrow Z_{j-1} \otimes S_1(F) \rightarrow \cdots
\rightarrow Z_1 \otimes S_{j-1}(F) \rightarrow S_j(F)
\rightarrow S_j(I) \rightarrow 0.
\]
be the graded pieces of the ${\mathcal Z}_{*}$-complex.
By Lemma~\ref{syzygetic}, it follows that $S_j(I) \simeq I^j$
for $1 \leq j \leq \ell-g+1$. Hence, from the above exact sequence resolving
$I^j$, we deduce, by induction on $j$, that ${\rm depth}\, Z_j \geq d-g+2$
for $1 \leq j \leq \ell-g$ as, by $(${\it b}$)$, ${\rm depth}\, I^j \geq
d-g-j+2$ for $1 \leq j  \leq \ell-g$.
Chasing depths in the defining exact sequences
\[
0 \rightarrow Z_{j+1} \rightarrow {\mathbb K}_{j+1} \rightarrow B_j
\rightarrow 0
\]
\[
0 \rightarrow B_j \rightarrow Z_j \rightarrow H_j(I)
\rightarrow 0
\]
we conclude that $H_j(I)$ has depth $d-g$, hence it is
Cohen--Macaulay, for $1 \leq j \leq \ell-g-1$. Clearly $H_0(I)
\simeq R/I$ is Cohen--Macaulay. By duality $($see \cite[2.13 and
2.22]{Hinv}$)$ this suffices to conclude that all Koszul homology
modules of $I$ are Cohen--Macaulay \QED

\medskip

In the following example we show that Lemma~\ref{syzygetic} and
Theorem~\ref{scm} break down if either the second analytic deviation
requirement or the depth conditions are relaxed.
The calculations have been carried out using the computer algebra system
{\tt Macaulay} \cite{BS}.
Example~\ref{failure}$(${\it b}$)$ was found by Johnson
\cite[2.12]{U} and has also played an important role in other
contexts $($e.g. \cite{Sri}$)$.

\begin{Example} \label{failure}
{\rm $(${\it a}$)$ \ The above results fail for ideals of higher
second analytic deviation even if the other assumptions are
satisfied: for instance, take $I$ to be the ideal generated by the
maximal minors of the matrix
\[
\left(
\begin{array}{cccc}
x_1 & x_2 & x_3 & x_4 \\
x_2 & x_3 & x_4 & x_1
\end{array}
\right)
\]
in the localized polynomial ring $R=k[x_1, \ldots, x_4]_{(x_1,
\ldots, x_4)}$.  In this case $I$ is a Cohen--Macaulay ideal with
$g=3$, $\ell = 4$ and $n=6$. However, $I$ is not even syzygetic, i.e.,
$S_2(I) \not\simeq I^2$.

\medskip

$(${\it b}$)$ \
The above results also fail if the depth conditions are not entirely
satisfied:
for instance, let $\varphi$ be an alternating $5 \times 5$ matrix of
variables and let ${\bf Y}$ be a $5 \times 1$ matrix of variables. Let $R$
be the localized polynomial ring obtained after adjoining all those 15
variables to a field $k$ of characteristic zero. The ideal $I=Pf_4(\varphi)
+I_1(\varphi \cdot {\bf Y})$, generated by the $4\times 4$ Pfaffians of
$\varphi$ and the entries of the product matrix $\varphi \cdot {\bf Y}$,
has the following properties
\begin{itemize}
\item[$($1$)$]
$I$ is a grade $5$ Gorenstein ideal with $\ell=9$ and $n=10$;

\item[$($2$)$]
$I$ has property $G_{\infty}$;

\item[$($3$)$]
$R/I^j$ is a Cohen--Macaulay ideal for $j=1, 2, 3$ while ${\rm depth}\,
R/I^4 =6$.
\end{itemize}
The ideal $I$ is, though, not even syzygetic.
For the above results to apply, ${\rm depth}\, R/I^4$ at least $7$
would be required. }
\end{Example}

\end{document}